\documentclass{ifacconf}

\usepackage{graphicx}      
\usepackage{natbib}        
\usepackage{cite}
\usepackage{graphicx,psfrag,epsfig,epsf}
\usepackage{amsmath}
\usepackage{amsfonts}
\usepackage{amssymb}
\usepackage{verbatim}
\usepackage{latexsym}
\usepackage{color}
\usepackage{subfigure}
\usepackage{algorithmic}
\usepackage{algorithm}
\usepackage{latexsym,theorem}
\usepackage{amsmath,amsfonts,amssymb,amsbsy}

\newcommand{\Fset}{\mathbb{F}}
\newcommand{\tr}{\top}
\newcommand{\Uset}{\mathbb{U}}

\newcommand{\Yset}{\mathbb{Y}}

\newcommand{\cS}{\mathcal{S}}

\newcommand{\cE}{\mathcal{E}}
\newcommand{\cN}{\mathcal{N}}

\newcommand{\tini}{\text{ini}}
\newcommand{\NLS}{\text{NLS}}
\newcommand{\ND}{\text{N-DeePC}}
\newcommand{\lin}{\text{lin}}
\newcommand{\bu}{\mathbf{u}}

\newcommand{\by}{\mathbf{y}}
\newcommand{\bg}{\mathbf{g}}
\newcommand{\bz}{\mathbf{z}}
\newcommand{\bb}{\mathbf{b}}
\newcommand{\bY}{\mathbf{Y}}

\newcommand{\bW}{\mathbf{W}}
\newcommand{\bU}{\mathbf{U}}

\newcommand{\col}{\operatorname{col}}

\begin{document}
\begin{frontmatter}

\title{Neural Data--Enabled Predictive Control}



\author[First]{Mircea Lazar}

\address[First]{Eindhoven University of Technology, Eindhoven, Netherlands \\ (e-mail: m.lazar@tue.nl)}

\begin{abstract}                
Data--enabled predictive control (DeePC) for linear systems utilizes data matrices of recorded trajectories to \emph{directly} predict new system trajectories, which is very appealing for real-life applications. In this paper we leverage the universal approximation properties of  neural networks (NNs) to develop neural DeePC algorithms for nonlinear systems. Firstly, we point out that the outputs of the last hidden layer of a deep NN  implicitly construct a basis in a so-called neural (feature) space, while the output linear layer performs affine interpolation in the neural space. As such, we can train off-line a deep NN using large data sets of trajectories to learn the neural basis and compute on-line a suitable affine interpolation using DeePC. Secondly, methods for guaranteeing consistency of neural DeePC and for reducing computational complexity are developed. Several neural DeePC formulations are illustrated on a  nonlinear pendulum example.
\end{abstract}

\begin{keyword}
Predictive control, Data--driven control, Nonlinear systems, Neural networks.
\end{keyword}

\end{frontmatter}
\section{Introduction}
\label{sec1}
Data--enabled predictive control (DeePC) \citep{Coulson2019} is a \emph{direct} data--driven control algorithm that utilizes data matrices to predict system trajectories and compute control actions without first identifying a parametric system model or predictor. This is very appealing for practical applications of model predictive control (MPC), as it simplifies the MPC design and implementation. For a detailed discussion on direct versus indirect data-driven methods for controller design we refer to \citep{Markovsky_2023}.

A unique feature of the DeePC algorithm is that it assigns to any given predicted/future input trajectory and past inputs and outputs a \emph{set} of predicted/future output trajectories, as explicitly pointed out in \citep[Equation (8)]{Lazar_2022}. This is in contrast with MPC or subspace predictive control (SPC) \citep{FAVOREEL_1999}, which assign a unique predicted output trajectory to every predicted input trajectory. As shown in \citep{Fiedler_2021}, in the case of noise free data, the set of trajectories spanned by the DeePC prediction equations collapses to the same unique trajectory as MPC and SPC. However, in the case of noisy data, when SPC yields an unbiased, least squares optimal output trajectory, the set of trajectories spanned by the DeePC prediction equations is richer, and it includes the SPC one. Hence, under suitable regularization \citep{Florian_Bridge}, DeePC can outperform MPC and SPC in the noisy data case, by optimizing the bias/variance trade--off. 

Currently, the design, implementation and stability analysis of linear DeePC are well understood, see, for example, the recent surveys \citep{Markovsky_2023, Handbook_2023} and the references therein. However, since most real-life applications are nonlinear, there is an increasing interest in developing DeePC formulations for nonlinear systems. Several promising results in this direction include using data--driven trajectory linearization techniques \citep{Berberich_2022}, kernel functions \citep{Huang_2024} or basis functions  expansion of input--output predictors \citep{lazar2023basis}. On the other hand, many \emph{indirect} data--driven predictive control (DPC) algorithms exist for nonlinear systems, which typically parameterize system dynamics or multi--step input--output predictors using deep neural networks (NNs), see, e.g., \citep{masti2020learning, Bonassi2021, Zarzycki2022, Lazar_CDC2023}. An advantage of using deep NNs is that one does not have to choose a specific kernel or basis function for representing the predictors, as the universal approximation theorem for deep NNs is valid for any Tauber-Wiener activation function \citep{Chen_1995_TW}. Also, there are many efficient toolboxes for training NNs \citep{Nelles_2020, Ljung_2020} and systematic guidelines for choosing the number of hidden neurons can be found in, e.g., \citep{Trenn_2008}. However, similarly to linear MPC or SPC, indirect neural DPC formulations assign to every future input trajectory a unique future output trajectory, corresponding to a nonlinear-least-squares-optimal neural predictor. Once trained off--line, such predictors are used on--line for prediction without exploiting newly measured data for improving predictions. Therefore, it would be of interest to develop neural-networks-based DeePC algorithms that exploit on--line measurements to optimize predictions. Alternatively, stochastic MPC methods with on-line NN model updates, see, e.g., \citep{Rolf_2023}, could be considered.

In this paper we therefore develop a neural DeePC algorithm for nonlinear systems by merging direct and indirect approaches to data-driven control \citep{Markovsky_2023}. Firstly, we point out that the outputs of the last hidden layer of a deep NN are implicitly forming a basis in a so-called \emph{neural (feature) space} (see, e.g., \citep{Widrow_2013}), while the output linear layer performs affine interpolation in the neural space. As such, we can train \emph{off-line} a deep NN using large data sets of trajectories to learn the neural basis and compute \emph{on-line} a suitable affine interpolation using DeePC. We show that for linear activation functions the original DeePC formulation is recovered, i.e., DeePC corresponds to a \emph{linear} NN mapping. 

The developed neural approach to DeePC has some desirable features compared to kernel-functions-based approaches  \citep{Huang_2024}, e.g., it decouples the dimension of the space in which on-line computations are performed from the number of data points. Moreover, compared to \citep{Huang_2024, lazar2023basis} it does not require choosing a special kernel or basis function (i.e., any Tauber-Wiener activation function can be chosen \citep{Chen_1995_TW}). Another contribution of this paper is the \emph{Neural}-DeePC-3 formulation defined in Problem~\ref{prob:DeePC3}, which provides a  very efficient regularization method for DeePC, i.e., auxiliary variables are restricted to $\Rset^{pN}$ (number of outputs times the prediction horizon) compared with regularization methods in \citep{Florian_Bridge, lazar2023basis} that employ auxiliary variables in $\Rset^T$ with $T>>pN$ the number of data points.

\section{Preliminaries}
\label{sec2}
Throughout this paper, for any finite number $q\in\Nset_{\geq 1}$ of vectors or functions $\{\xi_1,\ldots,\xi_q\}$ we will make use of the operator $\col(\xi_1,\ldots,\xi_q):=[\xi_1^\tr,\ldots,\xi_q^\tr]^\tr$. As the data generating system, we consider an unknown MIMO nonlinear system with inputs $u \in \Rset^m$ and measured outputs $y \in \Rset^p$. For example, such a system could be represented using a controllable and observable discrete-time state-space model:
\begin{equation}
 	\label{eq:2.1}
 	\begin{split}
 		x(k+1)&=\tilde f(x(k),u(k)), \quad k\in\Nset,\\
 		y(k)&=\tilde h(x(k)),\end{split}
 \end{equation}
where $x\in\Rset^n$ is an unknown state and $\tilde f$, $\tilde h$ are suitable, unkown functions. Given an initial condition $x(0)$ and a sequence of inputs $\bu_{[0,T-1]}:=\{u(0),\ldots,u(T-1)\}$, the system \eqref{eq:2.1} generates a corresponding output sequence $\by_{[1,T]}:=\{y(1),\ldots,y(T)\}$, which could be affected by measurement noise. 

\emph{Indirect} nonlinear data--driven predictive control, see, e.g., \citep{masti2020learning, Lazar_CDC2023}, typically uses a multi-step predictor of the NARX type, i.e., 
\begin{equation}
	\label{eq:2.2}
	\by_{[1,N]}(k):=\mathbb{F}(\bu_\tini(k),\by_\tini(k),\bu_{[0,N-1]}(k)),
\end{equation} 
where $\mathbb{F}:=\col(f_1,\ldots,f_N)$ and
\begin{align*}
	\bu_\tini(k)&:=\col(u(k-T_\tini),\ldots,u(k-1)) \in \Rset^{T_\tini m},\\ 
	\by_\tini(k)&:=\col(y(k-T_\tini+1),\ldots,y(k))\in \Rset^{
 T_\tini p},
\end{align*}
and where $T_\tini \in \mathbb{N}_{\geq 1}$ defines the order of the NARX dynamics (different orders can be used for inputs and outputs, but for simplicity we use a common order). Typically, the map $\Fset$ is parameterized using deep neural networks. Above, $\bu_{[0,N-1]}(k):=\{u(0|k),\ldots,u(N-1|k)\}$ and $\by_{[1,N]}(k):=\{y(1|k),\ldots,y(N|k)\}$ denote the sequence of predicted inputs and outputs at time $k$, respectively, based on measured data $\bu_\tini(k), \by_\tini(k)$. Note that since each $f_i$ is a MIMO predictor, it is the aggregation of several MISO predictors, i.e., $f_i=\col(f_{i,1},\ldots,f_{i,p})$ where each $f_{i,j}$ predicts the $j$-th output, i.e., for $i=1,\ldots,N$
\begin{equation}
\label{eq:MISO}
\begin{split}
	y_j(i|k)&=f_{i,j}(\bu_\tini(k),\by_\tini(k),\bu_{[0,N-1]}(k)),\\
	y(i|k)&=\col(y_1(i|k),\ldots,y_p(i|k)),
\end{split}
\end{equation} 
where $j=1,\ldots,p$ and $p$ is the number of outputs. 

For any $k\geq 0$ (starting time instant in the data vector) and  $j\geq 1$ (length of the data vector obtained from system \eqref{eq:2.1}), define
\begin{align*}
	\bar\bu(k,j) &:= \col(u(k),\ldots, u(k+j-1)),\\
    \bar\by(k,j)&: = \col(y(k), \ldots, y(k+j-1)).
\end{align*}
Then we can define the Hankel matrices:
\begin{equation}\label{eq:hankel_data}
	\begin{aligned}
		\bU_p &:= \begin{bmatrix}\bar\bu(0,T_\tini) & \ldots & \bar\bu(T-1, T_\tini) \end{bmatrix}, \\
		\bY_p &:= \begin{bmatrix}\bar\by(1, T_\tini) & \ldots & \bar\by(T, T_\tini) \end{bmatrix},\\
		\bU_f &:= \begin{bmatrix}\bar\bu(T_\tini, N) & \ldots & \bar\bu(T_\tini+T-1, N) \end{bmatrix},\\
		\bY_f &:= \begin{bmatrix}\bar\by(T_\tini+1, N) & \ldots & \bar\by(T_\tini+T, N) \end{bmatrix},
	\end{aligned}
\end{equation}
where $T\geq (m+p)T_\tini+mN$ is the number of columns of the Hankel matrices. 
\subsection{Deep neural networks}
\begin{figure}
	\centering
	\includegraphics[width=\columnwidth]{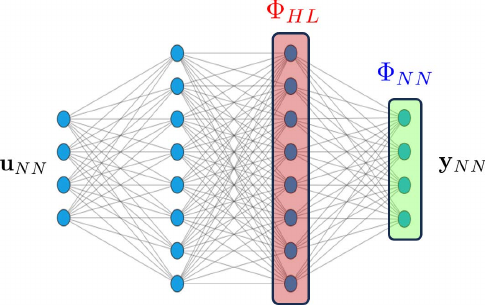}
	\caption{Illustration of a deep MLP NN.}
	\label{fig1}
\end{figure}

In this paper we will use deep multilayer perceptron (MLP)  neural networks \citep{Nelles_2020} for clarity of exposition, but in principle the adopted approach is applicable to other NN types as well. To define MLP NNs we consider arbitrary Tauber-Wiener nonlinear activation functions \citep{Chen_1995_TW} $\sigma:\Rset\rightarrow\Rset$ and linear activation functions, i.e., $\sigma_\lin(s):=s$. Given an arbitrary vector $x\in\Rset^n$ we define $\bar\sigma(x):=\col(\sigma(x_1),\ldots,\sigma(x_n))$, where $x_i$ are the elements of $x$. Let $\bu_{NN}$ and $\by_{NN}$ denote the NN inputs and outputs of suitable dimensions. Next, we can define a MLP NN map $\Phi_{NN}:\Rset^{q}\rightarrow \Rset^{p}$ recursively, i.e.,
\begin{equation}
\label{eq:NNmap}
\begin{split}
\bz_1&:=\bar\sigma(W_1\bu_{NN}+b_1),\\
\bz_i&:=\bar\sigma(W_i\bz_{i-1}+b_i),\,\, i=2,\ldots,l,\\
\by_{NN}&:=\bar\sigma_\lin(W_o\bz_{l}+b_o),
\end{split}
\end{equation}
where $\{(W_i,b_i)\}_{i=1,\ldots,l}$ are the matrices with weights and vectors with biases, respectively, of suitable dimensions, corresponding to each \emph{hidden} layer in the deep NN and $\{(W_o,b_o)\}$ are the matrix of weights and the vector of biases of the \emph{output} layer. The map $\Phi_{NN}(\bu_{NN})$ can then be defined as
\begin{equation}
\label{eq:NN}
\begin{split}
\Phi_{NN}(\bu_{NN})&:=\bar\sigma_\lin(W_o\bz_{l}+b_o)=W_o\bz_{l}+b_o\\
&=:W_o\Phi_{HL}(\bu_{NN})+b_o,
\end{split}
\end{equation}
where $\Phi_{HL}:\Rset^{q}\rightarrow \Rset^{L}$ is the map defined by the composition of hidden layers. $\Phi_{HL}$ maps the deep NN inputs from a space of dimension $q$ into a space of dimension $L$ (dictated by the number of neurons in the last hidden layer), which we refer to as the \emph{neural} space. Then the total deep NN map $\Phi_{NN}$ is obtained as an affine span in the neural space. See Figure~\ref{fig1} for an illustration. 

\begin{rem}
The definition \eqref{eq:NN} of $\Phi_{NN}$ as the affine span of the outputs of the hidden layer map $\Phi_{HL}$ shows that the hidden layers of a deep NN can be regarded as implicit basis functions generated by repeating two fundamental operations: affine span and passing through a scalar nonlinear activation function. In this way, by using a simple Tauber-Wiener activation function $\sigma$, a deep NN can generate a basis in the neural space, i.e., the space in which $\Phi_{HL}$ maps the NN inputs. Indeed, although not guaranteed in general, the outputs of the map $\Phi_{HL}$ are very likely to form a basis, i.e., they are linearly independent, as shown in, e.g., \citep{Chen_1995_TW, Widrow_2013}.  
\end{rem}

\section{Neural DeePC}
\label{sec3}
In order to build the neural DeePC controller, we need to define a neural network map $\Phi_{NN}:\Rset^{(m+p)T_\tini+mN}\rightarrow \Rset^{pN}$ that parameterizes the multi--step input--output NARX predictor $\Fset$ in \eqref{eq:2.2} and the corresponding hidden layer map $\Phi_{HL}:\Rset^{(m+p)T_\tini+mN}\rightarrow \Rset^L$, where $1\leq L\leq T$ is the number of neurons in the last hidden layer. Let $(\bW,\bb)$ denote all the weights and biases of all the hidden layers. Then given the Hankel data matrix
\[H:=\begin{bmatrix}\bU_p\\\bY_p\\ \bU_f \end{bmatrix}\in\Rset^{((m+p)T_\tini+mN)\times T},\]
let $H_{:j}$ denote the $j$-th column of $H$ and define the corresponding transformed matrix in the neural space as
\begin{equation}
\label{eq:3:data}
\bar\Phi_{HL}:=\begin{bmatrix}\Phi_{HL}(H_{:1}) & \Phi_{HL}(H_{:2}) &\ldots & \Phi_{HL}(H_{:T})\end{bmatrix}\in\Rset^{L\times T}.
\end{equation}
If $(\bW,\bb)$ are fixed, then $\bar\Phi_{HL}$ is a data matrix and, if $(\bW,\bb)$ are free variables, $\bar\Phi_{HL}$ is a matrix map/function. The last elements that need to be defined are the input and output vectors of the neural network, i.e.,
\[
\begin{split}
\bu_{NN}(k)&:=\col(\bu_\tini(k),\by_\tini(k),\bu_{[0,N-1]}(k))\\
\by_{NN}(k)&:=\by_{[1,N]}(k).
\end{split}
\]
Next, let $\Yset, \Uset$ be proper polytopic sets that represent constraints, let $l_s(y,u):=\|y-y_r\|_Q^2+\|u-u_r\|_R^2$ be a stage cost and $l_N(y)$ a terminal cost, taken for simplicity as $l_s(y,0)$, where $(y_r,u_r)$ are the desired output and input references and let $l_g(\bg(k))$ be a regularization cost. 

\begin{prob}[$Neural$-DeePC-1]
\label{prob:DeePC}
\begin{subequations}
	\label{eq:3.DeePC}
	\begin{align}
		\min_{\Xi(k)}\quad  & l_N(y(N|k))+\sum_{i=0}^{N-1} l_s(y(i|k), u(i|k))+\lambda l_g(\bg(k))\label{eq:3.DeePCa}\\
		&\text{subject to constraints:}\nonumber\\
		& \begin{bmatrix}\bar\Phi_{HL} \\ \bf{1}^\top\\\bY_f\end{bmatrix}\bg(k)=\begin{bmatrix}\Phi_{HL}(\bu_{NN}(k))\\1\\\by_{[1,N]}(k)\end{bmatrix}\label{eq:3.DeePCb}\\
		&(\by_{[1,N]}(k),\bu_{[0,N-1]}(k))\in\Yset^N\times\Uset^N\label{eq:3.DeePCc}
	\end{align}
\end{subequations}
\end{prob}
where $\Xi(k):=\col(\by_{[1,N]}(k),\bu_{[0,N-1]}(k),\bg(k))$ are the optimization variables and $\bf{1}$ is a vector of ones of the same length as $\bg(k)$, i.e., $T$. In principle, the hidden layers parameters $(\bW,\bb)$ can also be free variables that can be computed at every time $k$, but then it is difficult to guarantee that the resulting matrix $\begin{bmatrix}\bar\Phi_{HL} \\ \bf{1}^\top\end{bmatrix}$ will have full--row rank and Problem~\ref{prob:DeePC} becomes very complex, and most likely not solvable in real--time on--line. Therefore, we consider that $(\bW,\bb)$ are computed off-lline by solving the following nonlinear least squares (NLS) problem:
\begin{equation}
	\label{eq:3:NLS}
	(\bW^\ast,\bb^\ast,W_o^\ast,b_o^\ast):=
 \arg\min\left\|\bY_f-\bar\Phi_{NN}(\bW,\bb,W_o,b_o)\right\|_F^2
\end{equation}
where
$\bar\Phi_{NN}:=\begin{bmatrix}\Phi_{NN}(H_{:1}) & \Phi_{NN}(H_{:2}) &\ldots & \Phi_{NN}(H_{:T})\end{bmatrix}$.

Consider then the following least squares problem that can be used to recompute the weights and biases of the output layer after solving the NLS \eqref{eq:3:NLS}, i.e.
\begin{equation}
	\label{eq:3:LS}
	(W_o^{LS},b_o^{LS}):=
 \arg\min\left\|\bY_f-\begin{bmatrix}W_o&b_o\end{bmatrix}\begin{bmatrix}\bar\Phi_{HL}(\bW^\ast,\bb^\ast)\\\bf{1}^\top \end{bmatrix}\right\|_F^2.
\end{equation}
Notice that by solving the above least squares problem, i.e., after training the deep neural network using gradient descent and back propagation, we necessarily improve the data fit cost, since the new data fit cost will be either the same or lower. 

Next, we define the following multi--step input--output identified predictor:
\begin{equation}
\label{eq:3:pred}
\by_{[1,N]}^{\text{NLS}}(k):=\Phi_{NN}(\bu_{NN}(k),\bW^\ast,\bb^\ast,W_o^{LS},b_o^{LS}).
\end{equation}
As mentioned in the introduction, differently from the above predictor, the neural DeePC predictor is in general set--valued, i.e.
\begin{equation}
\label{eq:3:npred}
\begin{split}
&\by_{[1,N]}^{\text{N-DeePC}}(k)\in \\&\left\{\bY_f\bg(k) \ : \ \begin{bmatrix}\bar\Phi_{HL} \\ \bf{1}^\top\end{bmatrix}\bg(k)=\begin{bmatrix}\Phi_{HL}(\bu_{NN}(k))\\1\end{bmatrix}\right\},
\end{split}
\end{equation}
where the NLS optimal parameters $(\bW^\ast,\bb^\ast)$ are used in the map $\Phi_{HL}$. To analyze the relation between the above defined multi--step input--output predictors we introduce the following definition.
\begin{defn}
\label{def:eykoff}
Two models $\{M_1, M_2\}$ of system \eqref{eq:2.1} are called equivalent if for every constraints admissible input sequence $\bu_{[0,N-1]}$ and initial condition it holds that
\begin{equation}
\label{eq:echiv}
\|\by_{[1,N]}^{M_1}-\by_{[1,N]}\|=\|\by_{[1,N]}^{M_2}-\by_{[1,N]}\|,
\end{equation}
where $\by_{[1,N]}$ is the true system \eqref{eq:2.1} output. 
\end{defn}
Notice that one way to establish model equivalence is to show that $\by_{[1,N]}^{M_1}=\by_{[1,N]}^{M_2}$. 
Next, define
 \begin{align*}
 &\cS_\bg:=\\
 &\left\{\begin{bmatrix}\bar\Phi_{HL}\\\bf{1}^\top \end{bmatrix}^\dagger\begin{bmatrix}\Phi_{HL}(\bu_{NN}(k))\\ 1 \end{bmatrix}+\hat \bg\ : \ \hat \bg\in\cN\left(\begin{bmatrix}\bar\Phi_{HL}\\\bf{1}^\top \end{bmatrix}\right)\right\},
 \end{align*}
where $\cN(\cdot)$ denotes the null-space and $M^\dagger$ denotes the pseudo--inverse of $M$. 
\begin{thm}
\label{thm:eq}
Consider the nonlinear least squares optimal prediction model \eqref{eq:3:pred} and the neural DeePC prediction model \eqref{eq:3:npred} defined using the same set of data $\{\bU_p,\bY_p,\bU_f,\bY_f\}$ generated using system \eqref{eq:2.1} and the same map $\Phi_{HL}(\bW^\ast,\bb^\ast)$. Assume that the matrix $\begin{bmatrix}\bar\Phi_{HL}\\\bf{1}^\top \end{bmatrix}$ has full row--rank. Let $\cE:=\bY_f-\begin{bmatrix}W_o^{LS}&b_o^{LS}\end{bmatrix}\begin{bmatrix}\bar\Phi_{HL}\\\bf{1}^\top \end{bmatrix}$ be the matrix of residuals of the least squares problem \eqref{eq:3:LS}. Then the neural DeePC prediction model \eqref{eq:3:npred} is equivalent with the nonlinear least squares optimal prediction model \eqref{eq:3:pred} \emph{if and only if}  $\cE\hat\bg=\mathbf{0}$ for all $\hat \bg\in\cN\left(\begin{bmatrix}\bar\Phi_{HL}\\\bf{1}^\top \end{bmatrix}\right)$.
\end{thm}
\begin{pf}
From \eqref{eq:3:npred}, it follows that 
\[\begin{bmatrix}\bar\Phi_{HL} \\ \bf{1}^\top\end{bmatrix}\bg(k)=\begin{bmatrix}\Phi_{HL}(\bu_{NN}(k))\\1\end{bmatrix}\]
and thus all variables $\bg(k)$ that satisfy this system of equations satisfy $\bg(k)\in\cS_\bg$. Therefore, all predicted outputs generated by neural DeePC satisfy
\begin{align*}
&\by_{[1,N]}^{\text{N-DeePC}}(k)\in\\
&\{\bY_f\left(\begin{bmatrix}\bar\Phi_{HL}\\\bf{1}^\top \end{bmatrix}^\dagger\begin{bmatrix}\Phi_{HL}(\bu_{NN}(k))\\ 1 \end{bmatrix}+\hat \bg\right) :  \hat \bg\in\cN\left(\begin{bmatrix}\bar\Phi_{HL}\\\bf{1}^\top \end{bmatrix}\right)\}
\end{align*}
Then it holds that 
\begin{align*}
\bY_f\hat \bg&=\left(\cE+\begin{bmatrix}W_o^{LS}&b_o^{LS}\end{bmatrix}\begin{bmatrix}\bar\Phi_{HL}\\\bf{1}^\top \end{bmatrix}\right)\hat \bg\\&=\cE\hat \bg+\begin{bmatrix}W_o^{LS}&b_o^{LS}\end{bmatrix}\begin{bmatrix}\bar\Phi_{HL}\\\bf{1}^\top \end{bmatrix} \hat \bg=\cE\hat \bg.
\end{align*}
Since $(W_o^{LS},b_o^{LS})=\bY_f\begin{bmatrix}\bar\Phi_{HL}\\\bf{1}^\top \end{bmatrix}^\dagger$ it follows that
\[\by_{[1,N]}^{\text{N-DeePC}}(k)=\bY_f\begin{bmatrix}\bar\Phi_{HL}\\\bf{1}^\top \end{bmatrix}^\dagger\begin{bmatrix}\Phi_{HL}(\bu_{NN}(k))\\ 1 \end{bmatrix}=\by_{[1,N]}^{\text{NLS}}(k)\]
\emph{if and only if} $\cE\hat\bg=\mathbf{0}$ for all $\hat \bg\in\cN\left(\begin{bmatrix}\bar\Phi_{HL}\\\bf{1}^\top \end{bmatrix}\right)$.\qed
\end{pf}

In general the matrix of residuals $\mathcal{E}$ will not be identically zero, even in the case of noise free data, and thus, it is necessary to use a suitable regularization cost in neural DeePC, that penalizes the deviation of $\by_{[1,N]}^{\text{N-DeePC}}$ from $\by_{[1,N]}^{\text{NLS}}$. To this end, define the vector of variables \[\bg^{\NLS}(k):=\begin{bmatrix}\bar\Phi_{HL}\\\bf{1}^\top \end{bmatrix}^\dagger\begin{bmatrix}\Phi_{HL}(\bu_{NN}(k))\\ 1 \end{bmatrix}\]
and the regularization cost $l_g(\bg(k)):=\|\bg(k)-\bg^{\NLS}(k)\|_2^2$ in \eqref{eq:3.DeePCa}. Indeed, in this case for $\lambda\rightarrow\infty$ in \eqref{eq:3.DeePCa} we have that $\by_{[1,N]}^{\ND}(k)\rightarrow \by_{[1,N]}^{\NLS}(k)$. For a finite, sufficiently large value of $\lambda>0$ neural DeePC will optimize the bias/variance trade--off with respect to the  NLS optimal predictor \eqref{eq:3:pred}.

\begin{rem}It is relatively straightforward to observe that the original DeePC algorithm for linear (or affine) systems can be recovered as a special case of neural DeePC, i.e., Problem~\ref{prob:DeePC}. Indeed, by choosing $\sigma=\sigma_\lin$ for all the neurons in all the hidden layers, there exists a set of weights and biases $\{(W_i,b_i))\}_{i=1,\ldots,l}$ (the weights should be set equal to $0$ or $1$, and the biases should be set to $0$) such that $\Phi_{HL}(\bu_{NN})=\bu_{NN}$ and thus $\bar\Phi_{HL}=H$. Hence, the original DeePC algorithm corresponds to a \emph{linear} neural network map with the hidden layer parameters fixed to specific ($0$ or $1$) values. 
\end{rem}
\subsection{Computationally efficient neural DeePC formulations}
The regularization cost $l_g(\bg(k)):=\|\bg(k)-\bg^{\NLS}(k)\|_2^2$ is highly nonlinear and has a non--sparse structure, which leads to a computationally complex Problem~\ref{prob:DeePC}. Therefore, consider the alternative formulation of neural DeePC.

\begin{prob}[$Neural$-DeePC-2]
\label{prob:DeePC2}
\begin{subequations}
	\label{eq:3.DeePC2}
	\begin{align}
		\min_{\Xi(k)}\quad  & l_N(y(N|k))+\sum_{i=0}^{N-1} l_s(y(i|k), u(i|k))+\lambda l_g(\hat\bg(k))\label{eq:3.DeePC2a}\\
		&\text{subject to constraints:}\nonumber\\
		& \begin{bmatrix}\bar\Phi_{HL} \\ \bf{1}^\top\end{bmatrix}\hat\bg(k)=\begin{bmatrix}\bf{0}\\0\end{bmatrix}\label{eq:3.DeePC2b}\\
  	&\bY_f\left(\begin{bmatrix}\bar\Phi_{HL}\\\bf{1}^\top \end{bmatrix}^\dagger\begin{bmatrix}\Phi_{HL}(\bu_{NN}(k))\\ 1 \end{bmatrix}+\hat \bg(k)\right)=\by_{[1,N]}(k)\label{eq:3.DeePC2c}\\
  		&(\by_{[1,N]}(k),\bu_{[0,N-1]}(k))\in\Yset^N\times\Uset^N\label{eq:3.DeePC2d}
	\end{align}
\end{subequations}
\end{prob}
where $\Xi(k):=\col(\by_{[1,N]}(k),\bu_{[0,N-1]}(k),\hat\bg(k))$ are the optimization variables and $l_g(\hat\bg(k)):=\|\hat\bg(k)\|_2^2$. It is clear that the neural DeePC formulation in Problem~\ref{prob:DeePC2} is equivalent with the one in Problem~\ref{prob:DeePC}, as the predicted outputs $\by_{[1,N]}(k)$ are equivalently parameterized in the two problems. Moreover, when $\lambda\rightarrow\infty$ in Problem~\ref{prob:DeePC2} we obtain $\hat\bg(k)\rightarrow 0$ and thus $\by_{[1,N]}^{\ND}(k)\rightarrow \by_{[1,N]}^{\NLS}(k)$. However, in this case, the regularization cost is greatly simplified, which improves computational efficiency. 

Still, since  both $\bg(k)$ and $\hat\bg(k)$ are vectors  in $\Rset^{T}$, where $T$ is the data length that can be large, especially for nonlinear systems. Therefore we propose next an alternative formulation of Problem~\ref{prob:DeePC2}, where we aim to substitute $\bY_f\hat\bg(k)$ by a new vector of variables $\tilde\bg(k)\in\Rset^{pN}$. Note that if $\bY_f$ has full row--rank, we can parameterize $\hat\bg(k)$ as $\bY_f^\dagger\tilde\bg(k)$. We point out that this approach is novel also compared to the regularization methods proposed in \citep{lazar2023basis}.

\begin{prob}[$Neural$-DeePC-3]
\label{prob:DeePC3}
\begin{subequations}
	\label{eq:3.DeePC3}
	\begin{align}
		\min_{\Xi(k)}\quad  & l_N(y(N|k))+\sum_{i=0}^{N-1} l_s(y(i|k), u(i|k))+\lambda l_g(\tilde\bg(k))\label{eq:3.DeePC3a}\\
		&\text{subject to constraints:}\nonumber\\
		& \begin{bmatrix}\bar\Phi_{HL} \\ \bf{1}^\top\end{bmatrix}\bY_f^\dagger\tilde\bg(k)=\bf{0}\label{eq:3.DeePC3b}\\
  	&\bY_f\begin{bmatrix}\bar\Phi_{HL}\\\bf{1}^\top \end{bmatrix}^\dagger\begin{bmatrix}\Phi_{HL}(\bu_{NN}(k))\\ 1 \end{bmatrix}+\tilde \bg(k)=\by_{[1,N]}(k)\label{eq:3.DeePC3c}\\
  		&(\by_{[1,N]}(k),\bu_{[0,N-1]}(k))\in\Yset^N\times\Uset^N\label{eq:3.DeePC3d}
	\end{align}
\end{subequations}
\end{prob}
where $\Xi(k):=\col(\by_{[1,N]}(k),\bu_{[0,N-1]}(k),\tilde\bg(k))$ are the optimization variables and $l_g(\tilde\bg(k)):=\|\tilde\bg(k)\|_2^2$.

In Problem~\ref{prob:DeePC3} the free variables of DeePC are now mapped into the predicted output space of dimension $pN$, which is typically much smaller than $T$. The price to pay is less free variables for optimizing the bias/variance trade--off and it is necessary that the output data matrix $\bY_f$ has full row--rank. In the formulation of Problem~\ref{prob:DeePC3}, when $\lambda\rightarrow\infty$, we have that $\tilde\bg(k)\rightarrow 0$ and thus $\by_{[1,N]}^{\ND}(k)\rightarrow \by_{[1,N]}^{\NLS}(k)$. To avoid numerical instability, slack variables can be added on the right hand side of \eqref{eq:3.DeePC3b} and penalized in the cost function, as also done in linear DeePC.

\section{Illustrative example}
\label{sec4}
 Consider the following pendulum model \citep{Dpn2022}:
\begin{equation}
J \ddot{\theta}= u - b\dot{\theta}-\frac{1}{2}M L g \sin(\theta),
\end{equation}
where $u$ and $\theta$ are the system input torque and pendulum angle, while $J = \frac{ML^2}{3} kg\ m^2$, $M= 1 kg$ and $L = 1 m$ are the moments of inertia, mass and length of the pendulum. Moreover, $g = 9.81 m/s^2$ is the gravitational acceleration, and $b=0.1 Ns/m$ is the damping coefficient. One can model the pendulum dynamics using a state-space model discretized using Forward-Euler at $T_s = 0.033 $s as 

 \begin{equation}\label{eq:sys}
 \begin{split}
\begin{bmatrix}
x_1(k+1) \\
x_2(k+1) 
\end{bmatrix}&=\begin{bmatrix}
1-\frac{b T_s}{J} & 0 \\
T_s& 1
\end{bmatrix}  \begin{bmatrix}
x_1(k) \\
x_2(k) 
\end{bmatrix}+ \begin{bmatrix}
\frac{T_s}{J} \\
0 
\end{bmatrix} u(k)\\
&-\begin{bmatrix}
\frac{MLg T_s}{2J}\sin(x_2(k)) \\
0 
\end{bmatrix}, \\ 
    y(k)&=x_2(k),
 \end{split}    
 \end{equation}
where $x_1(k)$ is the angular velocity $\dot{\theta}$ and $x_2(k)$ is the angle $\theta$. We implemented and compared the performance and computational complexity of the 3 developed neural DeePC formulations, using the same data set, prediction horizon $N=10$ and tracking cost function weights $Q=200$ and $R=0.5$. For all regularization costs we have used $\lambda=1e+4=10^4$. For tuning the regularization weight parameter $\lambda$ we can make use of the Hanke-Raus heuristic \citep{Hanke96} and \citep[Proposition~II.1]{Lazar_2022}, which shows that linear DeePC with $\Pi$-regularization \citep{Florian_Bridge} can be equivalently written as a regularized least squares problem. The important difference with linear DeePC is that the control input sequence $\bu_{NN}(k)$ no longer depends linearly on $\bg(k)$ and as such, the input related stage cost function must be omitted in the tuning procedure. 

To generate the output data an open-loop identification experiment was performed using a multisine input constructed with the Matlab function \emph{idinput}, with the parameters \emph{Range} $[-4,\,4]$, \emph{Band} $[0,\,1]$, \emph{Period} $1000$, \emph{NumPeriod} $1$ and \emph{Sine} $[25,\,40,\,1]$. The data length is $1000$ and $T_\tini=5$ is used, as estimated in \citep{Dpn2022}. An MLP NN was defined in \emph{PyTorch} with one hidden layer with 30 neurons and with $tanh$ activation function. The network was trained to find the optimal weights and biases using the \emph{Adam} optimizer. The resulting map $\Phi_{HL}$ is used to generate a matrix $\bar\Phi_{HL}\in\Rset^{30\times 990}$, since $N=10$, the data length is $1000$ and we used  Hankel matrices to generate $(\bU_p,\bY_p,\bU_f)$. The resulting $\bar\Phi_{HL}$ matrix has full row--rank and a minimal singular value of $0.0188$.  

To compare the performance and computational complexity of all the derived data-driven predictive controllers we report the following performance indexes in Table~\ref{Tab:nonoise}: the integral squared error $J_{ISE}=\sum_{k=1}^{T_{\text{sim}}}\|y(k)-r(k)\|_2^2$, the integral absolute error $J_{IAE}=\sum_{k=1}^{T_{\text{sim}}}\|y(k)-r(k)\|_1$, the input cost $ J_u=\sum_{k=1}^{T_{\text{sim}}}\|u(k)\|_1$ and the tracking cost $J_{track}=\sum_{k=1}^{T_{sim}}\|Q^\frac{1}{2}(y(k)-y_r(k))\|_2^2+\|R^\frac{1}{2}(u(k)-u_r(k))\|_2^2$. The mean CPU time in seconds is also given in Table~\ref{Tab:nonoise}. All the predictive control optimization problem were solved with \emph{fmincon} and a standard laptop with Intel(R) Core(TM) i7-$8750$H CPU~$2.20$GHz, $16$GB RAM. 

\begin{table}[]
\centering
\begin{tabular}{|l|l|l|l|l|l|l|}
\hline
\textbf{Formulation} & $J_{ISE}$ & $J_{IAE}$  & $J_u$     & $J_{track}$  & CPU (s)\\ \hline
$Neural$-DeePC-1       &  16.24
     &     24.51     & 401 & 3759 & 91          \\ \hline
$Neural$-DeePC-2      &    16.32   &     24.62       &  401 &  3777 & 24       \\ \hline
$Neural$-DeePC-3     &16.32
    & 24.62        & 401 & 3776 & 0.39         \\ \hline
\end{tabular}
\caption{Performance \& mean CPU time.}
\label{Tab:nonoise}
\end{table}

\begin{figure}
	\centering
	\includegraphics[width=\columnwidth]{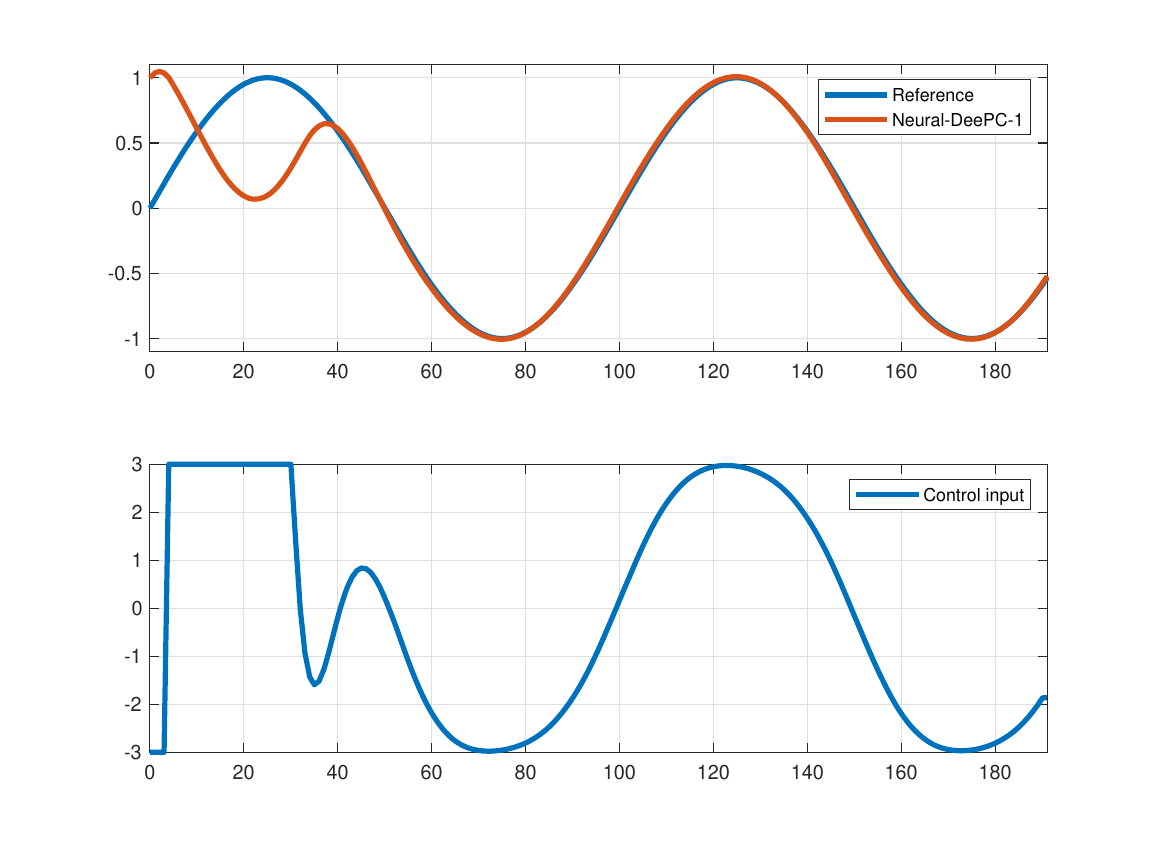}
	\caption{Tracking performance for $Neural$-DeePC-1.}
	\label{fig2}
\end{figure}

The closed--loop trajectories corresponding to $Neural$-DeePC-1 are plotted in Figure~\ref{fig2}.
\begin{figure}
	\centering
	\includegraphics[width=\columnwidth]{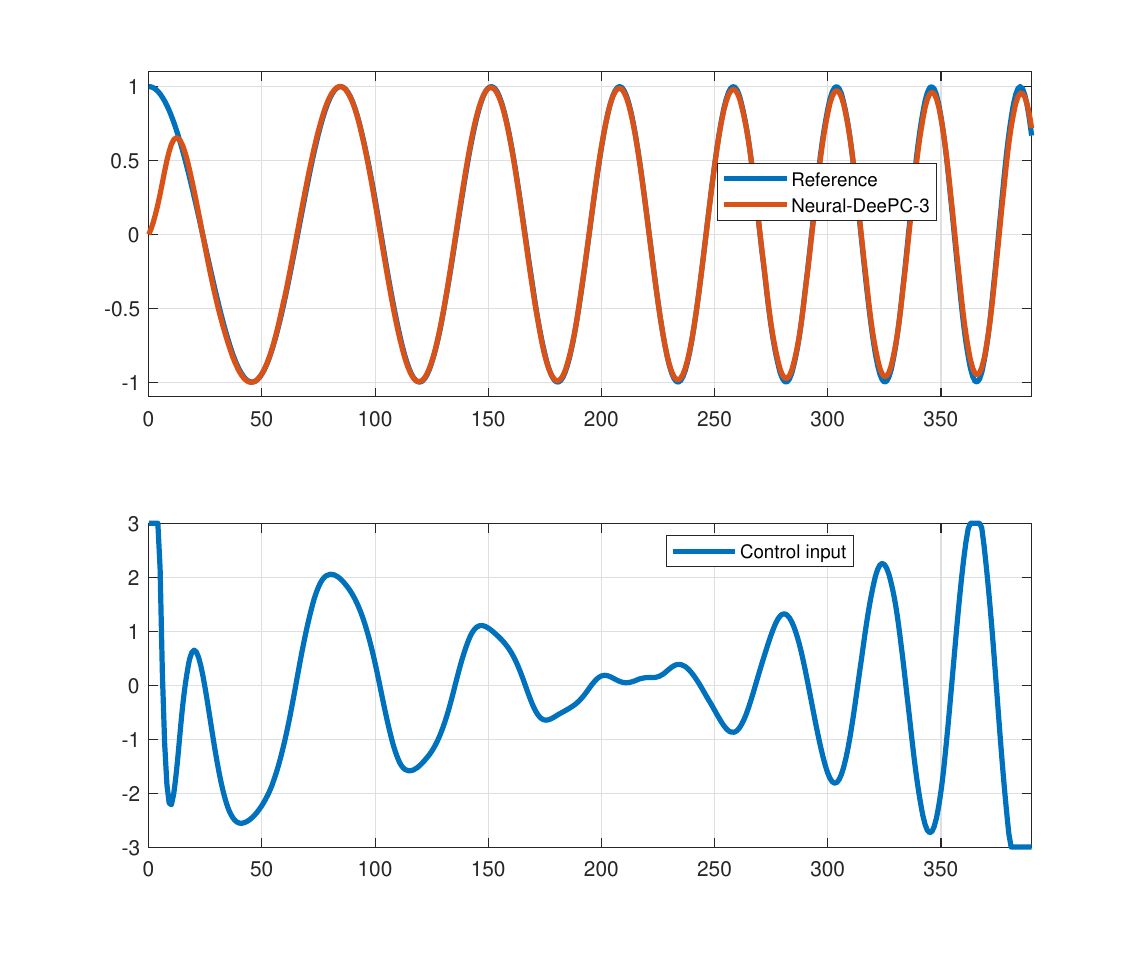}
	\caption{Tracking performance for $Neural$-DeePC-3.}
	\label{fig3}
\end{figure}
We observe that all the neural DeePC controllers perform well, since the obtained data is representative, the trained NN map is very accurate, and there is no noise. \emph{Neural}-DeePC-1 has slightly better performance, especially with respect to tracking performance. However, its computational complexity is much higher, which is expected, since it has $990$ additional optimization variables compared to an \emph{indirect} formulation and it uses a non-sparse regularization cost. Similar observations can be made for \emph{Neural}-DeePC-2, which ends up being about $4$ times faster than \emph{Neural}-DeePC-1. The \emph{Neural}-DeePC-3 formulation (implemented with slack variables in \eqref{eq:3.DeePC3b} with their squared 2-norm penalized in the cost function by $1e+4$) matches the performance of the other \emph{Neural}-DPC formulations, while achieving a much faster computation time, which is promising for real--life applications.

The performance of \emph{Neural}-DeePC-3 on a more challenging sinusoidal  reference with increasing frequency is shown in Figure~\ref{fig3}. In this case we have implemented Problem~\ref{prob:DeePC3} without any slack variables, which led to occasional warnings, but delivered a very promising average CPU time of $0.0748$ seconds.
\section{Conclusions}
\label{sec5}
In this paper we have constructed a formulation of data--enabled predictive control for nonlinear systems using deep neural networks. This was enabled by the key observation that the last hidden layer within a deep neural network implicitly constructs a basis for the neural network outputs. In the case of linear DeePC, the basis is simply formed by the system trajectories, which can be equivalently represented via a deep neural network with linear activation functions in the hidden layers. Future work will deal with data generation, training methods and architecture design for deep NNs with guarantees that the last hidden layer forms a representative basis in the space of output trajectories of dynamical systems.

\section*{Acknowledgements}
The author gratefully acknowledges the assistance of MSc. Mihai-Serban Popescu with constructing and training in PyTorch the deep NNs and with generating the corresponding map $\Phi_{HL}$ for the illustrative example. 

\bibstyle{ifacconf}
\bibliography{SysId_2024_ML}             








\end{document}